\input amstex\documentstyle{amsppt}  
\pagewidth{12.5cm}\pageheight{19cm}\magnification\magstep1
\topmatter
\title On the character of certain irreducible modular representations\endtitle
\author G. Lusztig\endauthor
\address{Department of Mathematics, M.I.T., Cambridge, MA 02139}\endaddress
\thanks{Supported in part by National Science Foundation grant DMS-1303060 and by a Simons Fellowship.}
\endthanks
\endtopmatter   
\document
\define\Rep{\text{\rm Rep}}

\define\Irr{\text{\rm Irr}}

\define\frl{\forall}

\define\sqc{\sqcup}

\define\qua{\quad}

\define\part{\partial}

\define\imp{\implies}
\define\ra{\rangle}
\define\n{\notin}
\define\iy{\infty}
\define\m{\mapsto}
\define\do{\dots}
\define\la{\langle}

\define\T{\times}

\redefine\i{^{-1}}

\define\Hom{\text{\rm Hom}}

\define\a{\alpha}

\redefine\d{\delta}
\define\e{\epsilon}

\define\r{\rho}

\redefine\l{\lambda}

\define\x{\xi}

\redefine\D{\Delta}

\define\kk{\bold k}

\define\pp{\bold p}
\define\qq{\bold q}

\define\NN{\bold N}

\define\ZZ{\bold Z}

\define\che{\check}

\define\AJS{AJS}
\define\FI{Fi}
\define\JA{Ja}
\define\KL{KL}
\define\SC{L1}
\define\MOD{L2}
\define\ST{St}
\define\WI{Wi}
\subhead 1\endsubhead
Let $G$ be an almost simple, simply connected algebraic group over $\kk$, an algebraically closed field of 
characteristic $p>1$. Let $\Rep G$ be the category of finite dimensional
$\kk$-vector spaces with a given rational linear action of $G$ and let $\Irr G$ be a set of representatives
for the simple objects of $\Rep G$. We fix a Borel subgroup $B$ of $G$ and a maximal torus $T$ of $B$; let
$Y=\Hom(\kk^*,T)$, $X=\Hom(T,\kk^*)$ (with group operation written as $+$) and let $\la,\ra:Y\T X@>>>\ZZ$ be
the obvious pairing. If $V\in\Irr G$ then there is a well defined $\l_V\in X$ with the following property: 
the $T$-action on the unique $B$-stable line in $V$ is through $\l_V$; according to Chevalley, $E\m\l_E$ is a
bijection from $\Irr G$ to a subset of $X$ of the form 
$$X^+=\{\l\in X;\la\che\a_i,\l\ra\in\NN\qua\frl i\in I\}$$ 
for a well defined basis $\{\che\a_i;i\in I\}$ of $Y$. For $\l\in X^+$ we shall denote by $V_\l$ the object 
of $\Irr G$ corresponding to $\l$.
If $V\in\Rep G$, then for any $\mu\in X$ we denote by $n_\mu(V)$ the multiplicity of $\mu$ in $V|_T$; we set
$[V]=\sum_{\mu\in X}n_\mu(V)e^\mu\in\ZZ[X]$ where $\ZZ[X]$ is the group ring of $X$ (the basis element of 
$\ZZ[X]$ corresponding to $\mu\in X$ is denoted by $e^\mu$ so that $e^\mu e^{\mu'}=e^{\mu+\mu'}$ for 
$\mu,\mu'\in X$). It is of considerable interest to compute explicitly the element $[V_\l]\in\ZZ[X]$ for any
$\l\in X^+$. Let $h$ be the Coxeter number of $G$.
A conjectural formula for $[V_\l]$ (assuming that $p\ge c_G^0$ where $c_G^0$ is a constant depending only
on the root datum of $G$) was stated in \cite{\SC, p.316}. 
In the early 1990's it was proved (see \cite{AJS} and the references there)
that there exists a (necessarily unique) prime number $c_G\ge c_G^0$ depending
only on the root datum of $G$ such that the conjectural formula in 
\cite{\SC, p.316} is true if $p\ge c_G$ and $c_G$ is minimum possible (but $c_G$ was not
explicitly determined). In \cite{\FI}, Fiebig showed that $c_G\le c'_G$ where
$c'_G$ is an explicitly known but very large constant. In \cite{\WI}, Williamson, 
partly in collaboration with Xuhua He, showed that for infinitely many $G$,
$c_G$ is much larger than $c_G^0$.
Now the conjecture in \cite{\SC, p.316} had an unsatisfactory aspect: it applied only 
to a 
finite set of $\l\in X^+$ which, after application of Jantzen's results \cite{\JA} on translation functors, 
becomes a larger but still finite set (including all $\l$ in 
$X^+_{red}=\{\l\in X^+;\la\che\a_i,\l\ra\le p-1\qua\frl i\in I\}$); then the case of 
a general $\l\in X^+$ had to be obtained by applying the Steinberg tensor product 
theorem \cite{\ST}. In this note I 
want to offer a reformulation of the conjecture in \cite{\SC, p.316} (now a known theorem for $p$ large 
enough) which applies directly to any $\l\in X^+$, see 7(b).

\subhead 2\endsubhead
{\it Notation.}
Let $NT$ be the normalizer of $T$ in $G$ and let $W=NT/T$ be the Weyl group. Note that $W$ acts naturally on
$T$ hence on $Y,X$ and $\ZZ[X]$. Let $\che R=\{y\in Y;y=w(\che\a_i)\text{ for some }w\in W,i\in I\}$ (the 
set of coroots). Define $\che\a_0\in\che R$ by the condition that $\che\a_0+\che\a_i\n\che R$ for any 
$i\in I$. (Thus $\che\a_0$ is the highest coroot). For $i\in I\sqc\{0\}$ define $\a_i\in X$ by the condition
that the map $X@>>>X$, $\l\m\l-\la\che\a_i,\l\ra\a_i$, is induced by a (uniquely defined) element $s_i$ of 
$W$. Let $w\m\e_w$ be the homomorphism $W@>>>\{1,-1\}$ such that $\e_{s_i}=-1$ for any $i\in I$. Define 
$\r\in X^+$ by $\la\che\a_i,\r\ra=1$ for any $i\in I$. Let $\le$ be the partial order on $X$ given by 
$\l\le\l'$ whenever $\l'-\l\in\sum_{i\in I}\NN\a_i$.

\subhead 3\endsubhead
Let
$$\D=\{\l\in X;\la\che\a_i,\l+\r\ra\le0\qua \frl i\in I,\la\che\a_0,\l+\r\ra\ge-p\}.$$
For $i\in I$ we define $s'_i:X@>>>X$ by $s'_i(\l)=\l-\la\che\a_i,\l+\r\ra\a_i$ (an affine reflection). We 
define $s'_0:X@>>>X$ by $s'_0(\l)=\l-(\la\che\a_0,\l+\r\ra+p)\a_0$ (an affine reflection). Let $W_a$ be the 
subgroup of the group of permutations of $X$ generated by $s'_i (i\in I\cup\{0\})$. Then $W_a$ is a Coxeter 
group on the generators $s'_i (i\in I\cup\{0\})$, with length function $l:W_a@>>>\NN$.

For $\l\in X$ we have $w\i(\l)\in\D$ for some $w\in W_a$ and among all such $w$ there is a unique one, 
$w_\l$, of minimal length. 

For $\l,\mu\in X^+$ we set $w=w_\l$ and
$$\pp_{\mu,\l}=\sum_{y\in W_a;y\i(\mu)=w\i(\l)}(-1)^{l(yw)}P_{y,w}(1)\in\ZZ$$
where $P_{y,w}$ is the polynomial associated in \cite{\KL} to $y,w$ in the Coxeter group $W_a$.

From the definitions we see that $\pp_{\mu,\l}\ne0\imp\mu\le\l$, $\pp_{\l,\l}=1$. Hence for $\l,\mu\in X^+$
we can define $\qq_{\mu,\l}\in\ZZ$ by the requirement
$$\sum_{\nu\in X^+}\pp_{\mu,\nu}\qq_{\nu,\l}=\d_{\mu,\l}$$
for any $\l,\mu$ in $X^+$. We have $\qq_{\mu,\l}\ne0\imp\mu\le\l$, $\qq_{\l,\l}=1$. 

\subhead 4\endsubhead
For any $\l\in X^+$ we can write uniquely $\l=\sum_{k\ge0}p^k\l^k$ where $\l^k\in X^+_{red}$ for all $k\ge0$
and $\l^k=0$ for large $k$. 

For any $\l\in X^+$ and any $k\in\NN$ we define elements $E^k_\l\in\ZZ[X]$ by induction on $k$ as follows:
$$E^0_\l=\sum_{w\in W}\e_we^{w(\l+\r)}/\sum_{w\in W}\e_we^{w(\r)},$$ 
$$E^k_\l=\sum_{\mu\in X^+}\pp_{\mu,\sum_{j;j\ge k-1}p^{j-k+1}\l^j}
E^{k-1}_{\sum_{j;0\le j\le k-2}p^j\l^j+p^{k-1}\mu}\text{ if }k\ge1.\tag a$$
Note that $E^k_\l\in\ZZ[X]^W$, the ring of $W$-invariants in $\ZZ[X]$.

We show that for $k\ge0$ we have
$$E^k_\l=\sum_{\mu\in X^+}\qq_{\mu,\sum_{j;j\ge k}p^{j-k}\l^j}
E^{k+1}_{\sum_{j;0\le j\le k-1}p^j\l^j+p^k\mu}.\tag b$$
For any $\mu\in X^+$ and any $k,h\ge0$,
$$(\sum_{j;0\le j\le k-1}p^j\l^j+p^k\mu)^h$$
is equal to $\l^h$ if $0\le h\le k-1$ and to $\mu^{h-k}$  if $h\ge k$; hence, by (a), we have
$$E^{k+1}_{\sum_{j;0\le j\le k-1}p^j\l^j+p^k\mu}=\sum_{\nu\in X^+}
\pp_{\nu,\sum_{h;h\ge k}p^{h-k}\mu^{h-k}}E^k_{\sum_{h;0\le h\le k-1}p^h\l^h+p^k\nu}.$$
Thus the right hand side of (b) is
$$\align&\sum_{\mu\in X^+}\qq_{\mu,\sum_{j;j\ge k}p^{j-k}\l^j}
\sum_{\nu\in X^+}\pp_{\nu,\sum_{h;h\ge k}p^{h-k}\mu^{h-k}}E^k_{\sum_{h;0\le h\le k-1}p^h\l^h+p^k\nu}\\&=
\sum_{\mu\in X^+}\qq_{\mu,\sum_{j;j\ge k}p^{j-k}\l^j}
\sum_{\nu\in X^+}\pp_{\nu,\mu}E^k_{\sum_{h;0\le h\le k-1}p^h\l^h+p^k\nu}\\&=
\sum_{\nu\in X^+}\d_{\nu,\sum_{j;j\ge k}p^{j-k}\l^j}E^k_{\sum_{h;0\le h\le k-1}p^h\l^h+p^k\nu}\\&=
E^k_{\sum_{h;0\le h\le k-1}p^h\l^h+p^k\sum_{j;j\ge k}p^{j-k}\l^j}\\&=E^k_\l.\endalign$$
This proves (b).

By induction on $k$ we see, using (a), that for any $k\ge1$ we have
$$\align&E^k_\l=\sum_{\mu_0,\mu_1,\do,\mu_{k-1}\text{ in }X^+}
\pp_{\mu_0,\l^0+p\mu_1}\pp_{\mu_1,\l^1+p\mu_2}\do\pp_{\mu_{k-2},\l^{k-2}+p\mu_{k-1}}\\&
\T \pp_{\mu_{k-1},\sum_{j\ge k-1}p^{j-k+1}\l^j}E_{\mu_0}^0.\tag c\endalign$$
By induction on $k$ we see, using (b), that for any $k\ge1$ we have
$$\align&E^0_\l=\sum_{\nu_0,\nu_1,\do,\nu_{k-1}\text{ in }X^+}
\qq_{\nu_0,\l}\qq_{\nu_1,(\nu_0-\nu_0^0)/p}\qq_{\nu_2,(\nu_1-\nu_1^0)/p}\do
\qq_{\nu_{k-1},(\nu_{k-2}-\nu_{k-2}^0)/p}\\&
\T E^k_{\nu_0^0+p\nu^0_1+\do+p^{k-2}\nu^0_{k-2}+p^{i-1}\nu_{k-1}}.\tag d\endalign$$

\subhead 5\endsubhead
Let $\l\in X^+$. We can find $n\ge0$ such that $\l^n=\l^{n+1}=\do=0$. If $k\ge n$ we have 
$$E^k_\l=E^{k+1}_\l.$$
Indeed, if $\mu\in X^+$ and $\qq_{\mu,\sum_{j;j\ge k}p^{j-k}\l^j}\ne0$ we have $\qq_{\mu,0}\ne0$ hence
$\mu\le0$, $\mu=0$ and $\qq_{0,0}=1$; it follows that 
$$E^k_\l=E^{k+1}_{\sum_{j;0\le j\le k-1}p^j\l^j}=E^{k+1}_\l.$$
Thus we can set $E^{\iy}_\l=E^k_\l$ for large $k$. Clearly, if $\l\in X^+_{red}$, then 
$E^1_\l=E^2_\l=\do=E^{\iy}_\l$. 

Letting $k\to\iy$ in 4(c),(d), we deduce that 
$$E^\iy_\l=\sum_{\mu_0,\mu_1,\mu_2,\do\text{ in }X^+;\mu_h=0\text{ for large }h}
(\pp_{\mu_0,\l^0+p\mu_1}\pp_{\mu_1,\l^1+p\mu_2}\pp_{\mu_2,\l^2+p\mu_3}\do)E^0_{\mu_0},$$
(note that for large $h$ we have $\pp_{\mu_h,\l^h+p\mu_{h+1}}=\pp_{0,0}=1$ so that the infinite product makes
sense) and 
$$E^0_\l=\sum_{\nu_0,\nu_1,\nu_2,\do\text{ in }X^+}(\qq_{\nu_0,\l}\qq_{\nu_1,(\nu_0-\nu_0^0)/p}
\qq_{\nu_2,(\nu_1-\nu_1^0)/p}\do)E^\iy_{\nu_0^0+p\nu^0_1+p^2\nu^0_2+\do}$$
(note that for large $h$ we have $\nu_h=0$ hence $\qq_{\nu_{h+1},(\nu_h-\nu_h^0)/p}=\qq_{0,0}=1$ so that the
infinite product makes sense).

\subhead 6\endsubhead
It is known since the early 1990's that, if $p$ is not very small, then the conjecture 8.2 in \cite{\MOD} on 
quantum groups at a $p$-th root of $1$ holds. In particular for $\l\in X^+$, the element $E^1_\l$ describes
the character of an irreducible finite dimensional representation of such a quantum group and the tensor
product theorem \cite{\MOD, 7.4} holds for it. 

Thus, if for any $\x=\sum_{\l\in X}c_\l e^\l\in\ZZ[X]$ (with $c_\l\in\ZZ$) and any $h\ge0$ we set
$\x^{(h)}=\sum_{\l\in X}c_\l e^{p^h\l}\in\ZZ[X]$, then for any $\l\in X^+$ we have the equality
$$E^1_\l=E^1_{\l^0}(E^0_{\sum_{j\ge1}p^{j-1}\l^j})^{(1)}.\tag a$$
We show by induction on $k\ge1$ that for any $\l\in X^+$ we have
$$E^k_\l=E^1_{\l^0}(E^1_{\l^1})^{(1)}\do(E^1_{\l^{k-1}})^{(k-1)}(E^0_{\sum_{j\ge k}p^{j-k}\l^j})^{(k)}.
\tag b$$
By (a), we can assume that $k\ge2$. Using 4(a), it is enough to show that
$$\align&\sum_{\mu\in X^+}\pp_{\mu,\sum_{j;j\ge k-1}p^{j-k+1}\l^j}
E^{k-1}_{\sum_{j;0\le j\le k-2}p^j\l^j+p^{k-1}\mu}\\&=
E^1_{\l^0}(E^1_{\l^1})^{(1)}\do(E^1_{\l^{k-1}})^{(k-1)}(E^0_{\sum_{j\ge k}p^{j-k}\l^j})^{(k)}.\endalign$$
Replacing here 
$$E^{k-1}_{\sum_{j;0\le j\le k-2}p^j\l^j+p^{k-1}\mu}=
E^1_{\l^0}(E^1_{\l^1})^{(1)}\do(E^1_{\l^{k-2}})^{(k-2)}(E^0_{\sum_{j\ge k-1}p^{j-k+1}\mu^j})^{(k-1)}$$
which is known from the induction hypothesis, we see that it is enough to show that
$$\align&\sum_{\mu\in X^+}\pp_{\mu,\sum_{j;j\ge k-1}p^{j-k+1}\l^j}
E^1_{\l^0}(E^1_{\l^1})^{(1)}\do(E^1_{\l^{k-2}})^{(k-2)}(E^0_{\sum_{j\ge k-1}p^{j-k+1}\mu^j})^{(k-1)}\\&
=E^1_{\l^0}(E^1_{\l^1})^{(1)}\do(E^1_{\l^{k-1}})^{(k-1)}(E^0_{\sum_{j\ge k}p^{j-k}\l^j})^{(k)}.\endalign$$
Thus, it is enough to show that
$$\sum_{\mu\in X^+}\pp_{\mu,\sum_{j;j\ge k-1}p^{j-k+1}\l^j}
(E^0_{\sum_{j\ge k-1}p^{j-k+1}\mu^j})^{(k-1)}
=(E^1_{\l^{k-1}})^{(k-1)}(E^0_{\sum_{j\ge k}p^{j-k}\l^j})^{(k)}$$
or that
$$\sum_{\mu\in X^+}\pp_{\mu,\sum_{j;j\ge k-1}p^{j-k+1}\l^j}E^0_{\sum_{j\ge k-1}p^{j-k+1}\mu^j})
=E^1_{\l^{k-1}}(E^0_{\sum_{j\ge k}p^{j-k}\l^j})^{(1)}.$$
Using (a), the right hand side is $E^1_{\sum_{j\ge k-1}p^{j-k+1}\l^j}$. This is equal to the left hand side,
by 4(a). This proves (b).

Letting $k\to\iy$ in (b) we obtain for any $\l\in X^+$:
$$E^\iy_\l=E^1_{\l^0}(E^1_{\l^1})^{(1)}(E^1_{\l^2})^{(2)}\do.\tag c$$
Note that for large $h$ we have $\l^h=0$ hence $E^1_{\l^h}=1$ so that the infinite product makes sense. (We
have also used that $E^0_{\sum_{j\ge k}p^{j-k}\l^j}=E^0_0=1$ for large $h$.)

\subhead 7\endsubhead
We now assume that $p\ge c_G$. Then, by the first paragraph of no.6, we have 
$$[V_\l]=E^1_\l\text{ for any }\l\in X^+_{red}.\tag a$$
Using the Steinberg tensor product theorem \cite{\ST} and (a), we see that for any $\l\in X^+$ we have
$$[V_\l]=[V_{\l^0}][V_{\l^1}]^{(1)}[V_{\l^2}]^{(2)}\do=E^1_{\l^0}(E^1_{\l^1})^{(1)}(E^1_{\l^2})^{(2)}\do.$$
Using this and 6(c) we deduce
$$[V_\l]=E^\iy_\l.\tag b$$

\subhead 8\endsubhead
We preserve the setup of no.7. In the identity
$$E^0_\l=\sum_{\mu\in X^+}\qq_{\mu,\l}E^1_{\mu}$$
(see 4(b)) the coefficient $\qq_{\mu,\l}$ can be interpreted as the multiplicity of an irreducible
representation of a quantum group at a $p$-th root of $1$ in a not necessarily irreducible representation of
that quantum group. In particular we have
$$\qq_{\mu,\l}\in\NN\tag a$$ 
for any $\l,\mu$ in $X^+$. We show by descending induction on $k$ that for any $\l\in X^+$ and any $k\ge0$
we have
$$E^k_\l=[V_\l(k)]\tag b$$
for some $V_\l(k)\in\Rep G$. 

If $k$ is large, we have $E^k_\l=E^\iy_\l$ and (b) follows from 7(b). Assume now that $k\ge0$ and that (b) 
is known when $k$ is replaced by $k+1$. Then (b) holds for $k$ by (a), 4(b) and the induction hypothesis.

\enddocument